\documentclass[amscd]{amsart}

\newtheorem{Thm}{Theorem}
\newtheorem{Cor}{Corollary}
\newtheorem{Lem}{Lemma}
\newtheorem{Prop}{Proposition}

\newtheorem{Fact}{Fact}

\theoremstyle{remark}
\newtheorem{Rem}{Remark}
\newtheorem{Def}{Definition}
\newtheorem{Ex}{Example}

\newsymbol\ltimes 226E
\newsymbol\rtimes 226F
\newsymbol\boxtimes 1202
\newsymbol\twoheadrightarrow 1310
\newsymbol\rightarrowtail 131A

\newcommand{\codim}{{\rm co}\dim}

\newcommand{\cal}{\mathcal}

\newcommand{\Hom}{{\rm Hom}}
\newcommand{\End}{{\rm End}}
\newcommand{\Ext}{{\rm Ext}}
\newcommand{\Ob}{{\rm Ob}}

\newcommand{\gr}{{\rm gr}}

\newcommand{\To}{\longrightarrow}
\newcommand{\iso}{{\widetilde \longrightarrow}}

\newcommand{\imbed}{\hookrightarrow}
\newcommand{\sur}{\twoheadrightarrow}

\newcommand{\bu}{\bullet}

\def\square{\hbox{\vrule\vbox{\hrule\phantom{o}\hrule}\vrule}}
\newcommand{\epf}{\square}

\newcommand{\bA}{{\bf A}}

\newcommand{\Ftil}{{\tilde {\mathcal F}}}

\newcommand{\N}{{\cal N}}
\newcommand{\Ntil}{{\tilde{\cal N}}}

\newcommand{\Nt}{{\tilde{\cal N}}}

\newcommand{\D}{{\mathfrak D}}

\newcommand{\DC}{{\mathbb DC}}

\newcommand{\oplusl}{\bigoplus\limits}
\newcommand{\cupl}{\bigcup\limits}

\renewcommand{\j}{{\frak j}}

\renewcommand{\L}{{\cal L}}
\renewcommand{\O}{{\cal O}}

\newcommand{\F}{{\cal F}}
\newcommand{\A}{{\cal A}}
\newcommand{\C}{{\cal C}}

\newcommand{\Ve}{{\cal V}}

\newcommand{\Fl}{{{\cal F}\ell}}

\newcommand{\GS}{{\cal S}}

\newcommand{\cons}{\underline{k}}

\newcommand{\nab}{\nabla}
\newcommand{\nabtil}{\tilde\nabla}
\newcommand{\del}{\Delta}
\newcommand{\MM}{M}
\newcommand{\NN}{N}
\newcommand{\sq}{\diamondsuit}

\renewcommand{\P}{{\cal P}}
\newcommand{\Ptil}{\tilde{\cal P}}

\newcommand{\Sh}{Sh}

\newcommand{\itil}{\tilde{i}}

\newcommand{\Obar}{\overline{O}}

\newcommand{\g}{{\mathfrak g}}
\newcommand{\gtil}{\tilde{\mathfrak g}}
\newcommand{\p}{{\mathfrak p}}

\newcommand{\Zet}{{\Bbb Z}}
\newcommand{\Ce}{{\Bbb C}}
\newcommand{\mult}{{{\Bbb C}^*}}

\newcommand{\Ibar}{\overline I}

\newcommand{\Db}{D^b}

\newcommand{\Gr}{{{\cal G}\frak r}} 

\newcommand{\pr}{{\rm pr}}

\newcommand{\<}{\langle}
\renewcommand{\>}{\rangle}

\newcommand{\nildat}{{\bf O}}

\newcommand{\bX}{{\bf X}}

\newcommand{\bY}{{\bf Y}}

 \newcommand{\I}{\iota}
 \newcommand{\Il}{\iota^l}
 \newcommand{\Ir}{\iota^r}
 
\newcommand{\Pil}{\Pi^l}
\newcommand{\Pir}{\Pi^r}

\newcommand{\proofnd}{{\it Proof }}

\title[Sheaves on nilpotent cone]{Quasi-exceptional sets and
equivariant coherent sheaves on the nilpotent cone}

\author{Roman Bezrukavnikov}

\begin{document}

\begin{abstract}
In \cite{izvrat} a certain $t$-structure on the derived category
of equivariant coherent sheaves on the nil-cone of a simple complex
algebraic group was introduced (the so-called perverse  $t$-structure
corresponding to the middle perversity). In the present note we show
that the same $t$-structure can be obtained from a natural quasi-exceptional
set generating this derived category. As a consequence we obtain a bijection
between the sets of dominant weights and pairs consisting of a nilpotent
orbit, and an irreducible representation of the centralizer of this
element, conjectured by Lusztig and Vogan (and obtained by other means
in \cite{B}).
\end{abstract}
\maketitle


\section{Introduction}

Let $G$ be a simple complex algebraic group,
 $\Lambda$ be the weight lattice of $G$,
and $\Lambda^+\subset \Lambda$ be the subset of dominant weights,
 $\g$ be the Lie algebra of $G$, and $\N\subset \g$ be the subvariety
of nilpotent elements.

 Let $\nildat$ be the set
of pairs $(O,L)$, where $O\subset  \N$ is a $G$-orbit, and $L$ is an
irreducible representation
of the centralizer $Z_G(n)$, $n\in O$ (up to conjugacy).

Lusztig and Vogan conjectured (independently) that there exists
 a natural bijection 
 between the sets $\nildat$ and $\Lambda^+$. (Since the
meaning of the word ``natural'' is not specified, this formulation
 of the conjecture is not precise).

Existence of such a bijection follows from the main result of \cite{B}
(the relation between the main result of \cite{B} and
the bijection $\nildat \leftrightarrow \Lambda^+$ is explained in
 \cite{cells4}, 10.8).
The argument of \cite{B} uses perverse sheaves on the affine flag
variety of the Langlands dual group, and some deep results of the geometric
theory of Langlands correspondence (in particular the construction of
\cite{KGB}). In this note we construct a bijection
$ \nildat \leftrightarrow \Lambda^+$
by more direct and elementary means.
(We do not check here
 that the bijection arising from the result of \cite{B}
coincides with the one constructed below; this will be done elsewhere).

Let us now describe the content of the paper.
 We provide a new (``exotic'')
$t$-structure on the triangulated categories $\Db(Coh^G(\N))$, the
derived category of $G$-equivariant coherent sheaves on $\N$. The
core $\P$ of this $t$-structure is an abelian category of finite
type (i.e. all objects have finite length); moreover it is a {\it
quasi-hereditary} (or {\it Kazhdan-Lusztig type}) category. This
means in particular, that $\P$ has a preferred ordered set of
objects called standard objects, another one of costandard
objects, and both these sets are in canonical bijection with the
set of (isomorphism classes of)
 irreducible objects (see section \ref{quaqua}
below for precise definitions).

The ``exotic'' $t$-structure admits two different descriptions.
On the one hand it is the perverse $t$-structure on equivariant coherent
 sheaves (in the sense of Deligne) corresponding to the middle perversity,
see \cite{izvrat}. This  makes clear  that
 (isomorphism classes of)
 irreducible objects in  $\P$ are numbered by  $\nildat$; for
 $(O,L)\in \nildat$ let $IC_{O,L}$ denote the corresponding irreducible
object of $\P$.

 On the other hand, given an abstract triangulated category $\D$ with
an ordered set of objects $\nab =\{ \nab^i \}$
satisfying certain conditions (a set satisfying those conditions
is called a dualizable quasi-exceptional
 set generating $\D$) one can produce a  $t$-structure  on $\D$, called the
 $t$-structure of a quasi-exceptional set. The core of this $t$-structure
is quasi-hereditary, and $\nab$ is the set of its costandard objects.
We show that for $\D=\Db(Coh^G(\N))$, the set $\nab=\{\nab^\lambda\}$,
$\lambda \in \Lambda^+$
consisting of direct images of positive line bundles under the Springer
map $\pi:\Ntil\to \N$ 
 is a quasi-exceptional set generating $\Db(Coh(N))$;
and that the corresponding $t$-structure coincides with the one described
in the previous paragraph.

 Thus the bijection
between the sets of irreducible and costandard objects in a quasi-hereditary
category yields a bijection
$$\nildat=\{IC_{O,L}\}\leftrightarrow \nab=\Lambda^+.$$

\medskip

Let us note that our approach is closely related to that of
\cite{O}. We recall briefly the set-up of {\it loc. cit.} Let
$K^G(\N)=K^0\left(Coh^G(\N)\right)$,  $K^{G\times \mult}(\N)
=K^0\left(Coh^{G\times \mult}(\N)\right)$ denote respectively the
Grothendieck groups of the category of $G$-equivariant coherent
sheaves on $\N$, and of the category of $G\times
\mult$-equivariant coherent sheaves on $\N$ (where $\mult$ acts on $\N$
by $t:n\to t^2n$).
 Then $K^G(\N)$  is freely generated by the classes of $\nab^\lambda$
($AJ_\lambda$ in notations of \cite{O}), $\lambda\in \Lambda^+$; and
$K^{G\times \mult}(\N)$ is freely generated as a
$\Zet[v,v^{-1}]$-module by the classes of $\nabtil_\lambda$. Here
$\nabtil_\lambda$ is a natural lift of $\nab^\lambda$ to
$Coh^{G\times \mult}(\N)$ (i.e. a $G\times \mult$ equivariant
coherent sheaf, which gives $\nab^\lambda$ upon restricting the
equivariance to $G$); and the action of the indeterminate $v$ on
$K^{G\times \mult}(\N)$ corresponds to the twist by the
tautological character of $\mult$. We call  the set of classes
$\{[\nabtil_\lambda]\}$ the costandard basis of $ K^{G\times
\mult}(\N)$.

 In \cite{O}
Ostrik conjectures existence of another $\Zet[v,v^{-1}]$-basis of
$K^{G\times \mult}(\N)$, which he calls the canonical basis. The
latter is characterized by properties similar to those
characterizing the Kazhdan-Lusztig basis in the (affine) Hecke
algebra. In particular  it is  in bijection with the standard
basis; and the transition matrix between the two bases is upper
triangular in a natural order, and is congruent to the identity matrix
 modulo $v^{-1}$.

In many known examples a $\Zet[v,v^{-1}]$ module with two bases as
above arises as a  Grothendieck group of a quasi-hereditary graded
category (see e.g. \cite{S}), with canonical and costandard
basis formed respectively by the classes of irreducible and
costandard objects, so one may ask whether this also happens in
the case under consideration. Indeed, a straightforward
generalization of our construction provides a $t$-structure on
$\Db(Coh^{G\times \mult}(\N))$ such that $\nabtil_\lambda$ are the
costandard objects of its core $\Ptil$, and $\tilde{IC}_{O,L}$ are
its irreducible objects;
here $\tilde{IC}_{O,L}$ is a natural lift of ${IC}_{O,L}$ to
$\Db(Coh^{G\times \mult}(\N))$. The classes of $\tilde{IC}_{O,L}$
form a basis of $K^{G\times \mult}(\N)$; this basis is obtained
from the costandard basis by an upper-triangular transformation.
In fact, $\{[\tilde{IC}_{O,L}]\}$ is the canonical
basis, whose existence is conjectured in \cite{O};
 however, this statement is not proved in this note, as I
was not able to find a  direct proof that $\Ptil$ is a 
{\it mixed category}
in the sense of \cite{BGS}, 4.1, which 
would guarantee that the
transformation matrix between the two bases is identity modulo
$v^{-1}$ (a known proof follows from the results of \cite{AB}).

\medskip

We remark that the methods of this note originate from the results
of \cite{AB}. In that paper we provide an equivalence between the
triangulated category $D^b(Coh^G(\N))$ and certain category
related to the derived category
$D^b(\Sh_\Fl)$, where  $\Sh_\Fl$ is the category of perverse sheaves on the
affine flag variety $\Fl$
of the Langlands dual group constant along the Schubert
stratification.
In particular,
we have an essentially surjective functor $F:D^b(\Sh_\Fl)\to D^b(Coh^G(\N))$.
Then $F$
sends the tautological
$t$-structure on  $D^b(\Sh_\Fl)$ (whose core is $\Sh_\Fl$)
 into the ``exotic''
$t$-structure on $D^b(Coh^G(\N))$ introduced in the
present paper.

Finally, we mention that results of \cite{AB} yield also
``exotic'' $t$-structures on the triangulated categories
$D^b(Coh^G(\Ntil))$, $D^b(Coh^G(St))$ (where $St=\Ntil\times_\N \Ntil$
 is the Steinberg
variety of triples). Those $t$-structures have Artinian (finite type)
 cores, and can
be described in terms of  (quasi)exceptional sets;
however, I do not know an analogue of the description  of the $t$-structures
 in terms of perverse
coherent sheaves, or precise structure of the
 irreducible objects of their cores.

\medskip

{\bf Acknowledgements.} I am much indebted to Leonid Positselskii
for inspiring discussions, and helpful critical comments.
I thank Victor Ostrik and George Lusztig for
their attention to the work; without 
their stimulating interest the paper may have never appeared.
I am also grateful to Bram Broer, Niels Lauritzen and Raphael Rouquier
for help with references. 

The author is supported by an NSF grant. 

\section{Quasi-exceptional sets and quasi-hereditary categories}\label{1}
Most of this section is a restatement of the result of \cite{BBD}
on glueing of $t$-structures. 
The results are most probably well known to the
experts, and appear in some form in the literature (cf e.g.
\cite{PS}; I have learnt many of them from
L.~Positselskii); we sketch the argument for the sake of
completeness. We work in a generality  slightly greater
 than usual (allowing possibly infinite exceptional sets), as this does
not require any additional efforts (for the application below
it would suffice to consider finite quasi-exceptional sets only).

\subsection{Quasi-hereditary categories}\label{quaqua}
 An abelian category $\A$
will be called of {\it finite type} if any object of $\A$ has
finite length. Let $\A$ be an abelian category of finite type
 with a fixed
ordering on the set $I$ of isomorphism classes of irreducible
objects.
 We fix a  representative $L_i$ in each isomorphism class
 $i\in I$.

For $n\in I$ let
 $\A_{\leq n}$, $\A_{<n}$ be the Serre subcategory in $\A$ generated
by  $L_i$ with $i\leq n$ or $i<n$ respectively.
 Thus $\A_{\leq n}$, $\A_{<n}$ are strictly full abelian
subcategories of
$\A$, and $X\in \A$ lies in $\A_{\leq n}$
(respectively in $\A_{<n}$) iff any irreducible subquotient
of $X$ is isomorphic to $L_i$ for some $i\leq n$ (respectively $i<n$).

\begin{Def} A pair $(\MM_n , \phi_n)$, where $\MM_n$ is an object
of $\A_{\leq n}$, and $\phi_n:\MM_n\to L_n$ is a nonzero
morphism is called a standard cover of $L_n$ if the following two
properties hold.

i) $Ker(\phi_n)\in \A_{< n}$.

ii) We have $\Hom(\MM_n, L_i)= 0=\Ext^1(\MM_n,L_i)$ for $i<n$.

\medskip

 A pair $(\NN^n , \phi^n)$, where $\NN^n\in \A_{\leq n}$,
 and $\phi^n:L_n\to \NN^n$ is a nonzero morphism is called
a costandard hull of $L_n$ if the following two properties hold.

i$'$) $CoKer(\phi^n)\in \A_{< n}$.

ii$'$) We have $\Hom(L_i,\NN^n)= 0=\Ext^1(L_i,\NN^n)$ for $i<n$.

We will say that an object $M$ is standard (costandard)
if some morphism to (from) an irreducible object from (to) $M$
is a standard cover  (respectively, a costandard hull).
\end{Def}

\begin{Lem} A (co)standard cover (hull) is unique up to a unique isomorphism if
 exists.
\end{Lem}

\proof If $\phi_n:\MM_n \to L_n$ and $\phi_n':\MM_n'\to L_n$ are
two standard covers. Then $\Hom(\MM_n,
Ker(\phi_n'))=0=Ext^1(\MM_n, Ker(\phi_n'))$ implies that
$$\Hom(L_n, \MM_n')\iso \Hom(\MM_n,\MM_n').$$ In particular there
exists a unique morphism $\MM_n\to \MM_n'$ compatible with
$\phi_n$, $\phi_n'$, which proves the claim about standard covers.
The argument for costandard hulls is parallel. \epf

\begin{Def} A quasi-hereditary category is a finite type abelian category
with an ordering on the set of isomorphism classes of irreducible
objects, such that a standard cover, and an irreducible hull exist
for any irreducible object of $\A$.
\end{Def}

\subsection{Quasi-exceptional sets}\label{Qes}
We first fix some notations, partly borrowed from \cite{BBD}.
 Let $\D$ be a triangulated category.
We will write $\Hom^n(X,Y)=\Hom(X,Y[n])$, and  denote the
graded abelian group $\oplusl_n \Hom(X, Y[n])$ by   $\Hom^\bu(X,Y)
$;
 also

\noindent $\Hom^{>0}(X,Y)=\oplusl_{n>0} \Hom(X,Y[n])$
etc.

For an object $X$ of a category we will write $[X]$ for its
isomorphism class. For a category $\C$ let $[\C]$ be the set
of isomorphism classes of $\C$.

Let $\D$ be a triangulated category.

 If $\bX,\bY$ are
subsets of  $[\D]$, then $\bX*\bY$ denotes the subset of
$[\D]$ consisting of
 classes of all objects $Z$, for which there exists an
exact triangle $X\to Z\to Y\to X[1]$ with $[X]\in \bX$, $[Y]\in
\bY$. The octahedron axiom implies (see \cite{BBD},
 Lemma 1.3.10) that the $*$-operation is associative, so $\bX_1*\bX_2\cdots
*\bX_n$ makes sense. For a subset $\bX\subset [\D]$
 let $\<X\>\subset \D$ be the strictly
full subcategory defined by
$$[\< \bX\>]= \cupl_n \bX * \bX *\cdots * \bX,$$
where $\bX$ appears $n$ times in the right-hand side.


For $\bX\subset [\D]$ the triangulated subcategory generated
by $\bX$ is the smallest strictly full triangulated subcategory
$\D_{\bX}\subset \D$, such that $[\D_{\bX}]\supset \bX$.
 Thus
$$\D_\bX=\< \cupl_{n\in \Zet} \bX[n]\>,$$
where $\bX[n]=\{\left[ X[n] \right] \ ,\ [X]\in \bX\}$. We will say
``the triangulated subcategory generated by objects/subcategories''
instead of ``the triangulated category generated by the corresponing set
of isomorphism classes'', and write $\<X\>$ instead of $\< \{[X]\} \>$ etc.

\medskip

An ordered subset  $\nab=\{ \nab^i,\;  i\in I\}$ of $\Ob(\D)$ is
called {\it quasi-exceptional} if we have $\Hom^\bu (\nab^i,
\nab^j)=0$ for $i<j$,
 $\Hom ^{<0}(\nab^i,\nab^i)=0$, and $\End(\nab^i)$ is a division algebra
for all $i$.


For a quasi-exceptional set $\nab$, and $\nab^i\in \nab$ we set
  $\D_{\leq i}=\D_{\{\nab _j\ |\ j\leq i \} }$, $\D_{<i}=\D_
{\{\nab _j\ |\ j < i\}}$.

For a full triangulated subcategory $\D'\subset \D$ we will denote
 by $\D/\D'$ the factor category;
then $\D/\D'$ is again a
triangulated category (see \cite{V}, 2.2.10). For $X,Y\in \D$ we will
denote by $X\mod \D'$ the image of $X$ in $\D/\D'$,
and will write $X\cong Y\mod \D'$ instead of
$(X\mod \D')\cong (Y\mod \D').$

Let $\nab= \{ \nab^i,\; i\in I\}$ be a quasi-exceptional set, and
let  $\del=\{ \del_i,\;  i\in I\}$ be another subset of
$\Ob(\D)$ (in bijection with $\nab$).

We say that $\del$ is  {\it dual} to $\nab$ if
 \begin{equation}\label{YX}
\Hom^\bu (\del_n, \nab^i)=0 \ \ {\rm for}\ \ n>i;
\end{equation}
and there exists an isomorphism
 \begin{equation}\label{XY}
\del_n \cong \nab^n \mod \D_{<n}.
\end{equation}
We set  $$\sq_n \overset{\rm def}{=}
\del_n\mod \D_{<n} \cong \nab^n \mod \D_{<n}
\in \D_{\leq n}/\D_{<n}$$

\begin{Lem}\label{dual}
If $\nab$ is a
 quasi-exceptional set, and $\del$ is a dual set,
then

a)  $\Hom^\bu(\del_i,X)=0=\Hom^\bu(X,\nab^i)$ for all $X\in \D_{<i}$.

b) $\Hom^\bu (\del_i, \nab^j)=0$ unless $i=j$.

c) For all $X\in \D_{\leq i}$ we have
\begin{equation}\label{homquot}
\begin{array}{ll}
\Hom^\bu(\del_i, X)\iso \Hom^\bu(\sq_i, X\mod \D_{<i});\\
\Hom^\bu(X,\nab^i)\iso \Hom^\bu(X\mod \D_{<i}, \sq_i).
\end{array}
\end{equation}

d) $\Hom^\bu(\del_i,\del_j)=0$ for
$i>j$, and
 \begin{equation}\label{delnab}
 \Hom^\bu(\del_i,\del_i)\cong \Hom^\bu(\sq_i,\sq_i)\cong
\Hom^\bu (\nab^i,\nab^i)\cong \Hom^\bu(\del_i,\nab^i).
\end{equation}
 The induced isomorphisms
 $\End(\del_i)\cong \End(\nab^i)\cong \End(\sq_i)$ are isomorphisms
of algebras.

e) Let  $\nab$ be a quasi-exceptional set, and $\del$, $\del'$ be
two dual sets. Then
\begin{equation}\label{deldel}
\del_i\cong \del'_i \end{equation}
 for all $i$; moreover, there
exists a unique isomorphism \eqref{deldel} compatible with a fixed
 isomorphism \eqref{XY}.

f) Assume $\nab$ is well-ordered (i.e. every subset of $\del$
has the minimal element).
Then we have
 $\D_{\leq n}=\D_{\{\del_i| i\leq n\} }$,
 $\D_{< n}=\D_{\{\del_i| i< n\} }$. In particular,
if $\nab$  generates $\D$ (as a triangulated category), then so does
 $\del$.

\end{Lem}

\proof (a) is immediate from the definition.

If $i>j$, then (b) follows from the first equality in (a); while
if $i<j$, then it follows from the second equality in (a).

By
 \cite{V}, chapitre II, Proposition 2.3.3(a) part (a) of the Lemma implies
part (c).

 (d) and (e) follow from (c).

Finally (f) follows from the definition by (transfinite) induction. \epf

 \begin{Rem}
Let  $\nab= \{ \nab^i|\;  i\in I\}$ be a quasi-exceptional set,
and let $\Ibar $ be the set $I$ with the opposite ordering.
Statement (d) of the Lemma shows that  if $\del$ is a dual set
for   $\nab$, then $\del$ is a quasi-exceptional set indexed by
$\Ibar$.
 We say that  a quasi-exceptional set $\nab$ is {\it dualizable},
 if a dual set  exists.
\end{Rem}

 \begin{Rem} A quasi-exceptional set is called {\it exceptional}
(see e.g. \cite{BK})
if $\D$ is $k$-linear for a field $k$, and $\Hom^{\bu}(\nab^i,\nab^i)=k$.
 It is proved in \cite{BK} that if $\Hom^\bu(X,Y)$ is a finite dimensional
$k$ vector space for all $X,Y\in \D$, then
any (finite) exceptional set in $\D$ is dualizable.
\end{Rem}

\begin{Ex}\label{prevrat}
The reader can keep in mind the following example.
Let $\D$ be a full subcategory in the bounded  derived category of
 sheaves of $k$-vector spaces
 on a reasonable topological space (or of etale sheaves on a reasonable
scheme),
consisting of complexes whose cohomology is
smooth along a fixed stratification. Assume for simplicity that
the stratum $\Sigma_i$ are connected and simply-connected; we
write
$j<i$ if $\Sigma_j$ lies in the closure of $\Sigma_i$.
 Let
$\j_i$  denote the imbedding of  $\Sigma_i$ in the space.
Let $p_i$ be
arbitrary integers.
Then  objects $\nab^i=\j_*(\cons [p_i])$ form a quasi-exceptional set
generating
$\D$, and $\del_i=\j_!(\cons [p_i])$ is the dual set.
$\nab,\del$ are exceptional iff the strata are acyclic.
\end{Ex}

\begin{Prop}\label{tstr} Let $\D$ be a triangulated category.
Let $I$ be a totally ordered set, and $\nab^i$, $i\in I$ be a
dualizable quasi-exceptional set in $\D$, which generates $\D$ as a
triangulated category;
let $\del$   be the dual set.

 There exists a unique $t$-structure $(\D^{\geq 0}, \D^{<0})$ on
$\D$, such that $\nab^i\in \D^{\geq 0}$; $\del^i\in \D^{\leq 0}$.
Moreover, $\D^{\geq 0}$, $\D^{<0}$ are given by
\begin{equation}\label{defDbo}
  \D^{\geq 0}=\< \{\nab^i[d]\ ,\ i\in I, \; d\leq 0\} \>;
\end{equation}
\begin{equation}\label{defDme}
\D^{< 0}=\< \{\del_i[d]\ ,\ i\in I, \; d> 0\} \>.
\end{equation}

\end{Prop}

We will need two Lemmas to prove the Proposition. The first one
settles the case when $\nab=\del$ consists of one element (``the base of
induction''); the second one allows to use the gluing of $t$-structure
(see \cite{BBD}) to
make an induction step.

\begin{Lem}\label{imposic} \cite{P}
a) Let $\D$ be a triangulated category, and $\A\subset \D$ be a full semisimple
abelian
 subcategory, which generates $\D$ as a triangulated category.
Suppose that
\begin{equation}\label{nootr}
\Hom^{<0}(X,Y)=0\ \ \ {\rm for} \ \ X,Y\in \A.
\end{equation}
Then there exists a unique $t$-structure on $\D$ whose core contains $\A$;
it is given by:
\begin{equation}\label{ssA}
\begin{array}{lll}
\D^{\leq 0}=\< \cupl_{ i \geq 0} \A[i]\>;\\
\D^{\geq 0}=\< \cupl_{ i \leq 0} \A[i]\>;\\
\D^{\leq 0}\cap \D^{\geq 0}=\< \A\>.
\end{array}
\end{equation}

b) The set of isomorphism classes of simple objects of $\A$
coincides with the set of isomorphism classes of simple objects of
$\<\A\>$.
\end{Lem}

\proof
(a) follows from \cite{BBD}, Remarque 1.3.14.
More precisely, {\it loc. cit.} shows that the conclusion of (a)
holds for any
 full subcategory $\A\subset \D$ which satisfies \eqref{nootr} and such that
\begin{equation}\label{admiBBD}
[\A]*[\A][1] \subset [\A][1] *[\A]
\end{equation}
 (a subcategory satisfying \eqref{nootr}, \eqref{admiBBD}
is called {\it admissible} in \cite{BBD}, Definition 1.2.5).
 A semisimple full
abelian subcategory satisfying \eqref{nootr} is readily seen to be
admissible in this sense; indeed, for such a category we have
$$ [\A]*\left[\A[1]\right]=\{X\oplus Y[1]\ |\ X,Y\in \A\}.$$

Recall from \cite{BBD}, Proposition 1.2.4 that a sequences $0\to
X\overset{f}{\to} Y\overset{g}{\to} Z\to 0$ in an admissible
abelian subcategory in $\D$ is exact iff there exists a
distinguished triangle
$$X\overset{f}{\to} Y\overset{g}{\to} Z\overset{h}{\to} X[1];$$
in particular this is true for the subcategory $\<\A\>=\D^{\leq
0}\cap \D^{\geq 0}$, as the core of any $t$-structure is
admissible. Hence
  every object of $\<\A\>$ has a finite
filtration whose subquotient are simple in $\A$.
It remains to see that these objects are also simple in $\<\A\>$.
But if $L\in \A$ is not simple in $\<\A\>$,
 then there exists a simple object $L'\in \A$, and a nonzero morphism
$L'\to L$ which is not an isomorphism; so $L$ is not simple in $\A$.
\epf

\begin{Cor}\label{el1}
If $\D=\D_{\{X\}}$ for an object $X\in \D$ such that
$\Hom^{<0}(X,X)=0$, and $\End(X)$ is a division algebra, then
there exists a unique $t$-structure on $\D$ whose core contains $X$.
It is given by
\begin{equation}\label{adynelement}
\begin{array}{lll}
\D^{\leq 0}=\< \{  X[i]\>\ |\  i \geq 0\} ;\\
\D^{\geq 0}=\< \{  X[i]\>\ |\  i \leq 0 \} \>;\\
\D^{\leq 0}\cap \D^{\geq 0}=\< X\>.
\end{array}
\end{equation}
$X$ is a simple object of the core of this $t$-structure.
\end{Cor}

\proof Apply the previous Lemma to
 $\A=\{ X^{\oplus n}\ |\ n\in \Zet_{\geq 0}\}$. \epf

For a subcategory
$\A$ in an additive category $\C$ let us (following \cite{BK}) write
 $\A^\perp=(\A^\perp)_\C$ (respectively $^\perp \A=(^\perp \A)_\C$)
for a strictly full subcategory in $\C$ consisting of objects $X$ for which
$\Hom(A,X)=0$ (respectively $\Hom(X,A)=0$) for all $A\in \A$. The subcategories
$\A^\perp$, $^\perp \A$ are called respectively  right and left orthogonal
of $\A$.

Set $\D_n=\D_{\{\del_n\} }$, $\D^n=\D_{\{ \nab^n\} }$.

\begin{Lem}\label{2lem} a) We have
$$ [\D_{\leq n}]=[\D_n]*[\D_{<n}];$$
$$
[\D_{\leq n}]=[\D_{<n}] * [\D^n].
$$

b) We have $$\D_n=(^\perp \D_{<n})_{\D_{\leq n}};$$
 $$\D^n=( \D_{<n}^\perp)_{\D_{\leq n}};$$
$$(\D_n^\perp)_{\D_{\leq n}}=\D_{<n}= (\D_n^\perp)_{\D_{\leq n}}.$$

c) $\D_{<n}$ is a thick (saturated) subcategory in $\D_{\leq n}$.

d) The projection
 $\Pi=\Pi_n: \D_{\leq n}\to \D_{\leq n}/\D_{<n}$
induces equivalences of triangulated categories
$$\Pi|_{\D_n}: \D_n\iso \D_{\leq n}/\D_{<n};$$
$$\Pi|_{\D^n}:\D^n\iso \D_{\leq n}/\D_{<n}.$$
$\Pi$ has a left adjoint $\Pil$ and a right adjoint $\Pir$.
Moreover, $\Pil$ maps $\D_{\leq n}/\D_{<n}$ to $\D_n$ and induces
an equivalence inverse to $\Pi|_{\D_n}$; while $\Pir$
 maps $\D_{\leq n}/\D_{<n}$ to $\D^n$ and induces
 an equivalence inverse to $\Pi|_{\D^n}$.

e) The inclusion functor $\I:\D_{<n}\imbed \D_{\leq n}$ has a left adjoint
$\Il$ and a right adjoint $\Ir$. Functors $\Il$, $\Ir$ are triangulated
(i.e. send distinguished triangles into distinguished triangles).
\end{Lem}

\proof
It is obvious that $\D_n\subset (^\perp \D_{<n})_{\D_{\leq n}}$,
 $\D^n=( \D_{<n}^\perp)_{\D_{\leq n}}$. Hence
$$[\D_{<n} ]* [\D_n] \subset [\D_n]*[\D_{<n}];$$
$$[\D^n]*[\D_{<n}]\subset [\D_{<n}] * [\D^n].$$
Now (a) follows from the fact that  $\D_{\leq n}$ is generated
by $\D_n$, $\D_{<n}$, as well as by $\D_{<n}$, $\D^n$, by associativity of
the star operation.

(b) is immediate from (a); (c) follows from (b) because both left
and right orthogonal of a triangulated category is a thick
subcategory. The rest of the Lemma follows e.g. from \cite{V},
chapitre II, Proposition 2.3.3 (see also \cite{BK} 1.5--1.9).

Let us recall the construction of adjoint functors $\Pil$, $\Pir$,
$\Il$, $\Ir$ for further reference. By part (a) of the Lemma for $X\in
\D_{\leq n}$ there
 exist distinguished
triangles $$X_n \to X \to X_{<n}\to X_n[1];$$
 $$X^{<n} \to X \to X^{n}\to X^{<n}[1]$$
with $X_n\in \D_n$, $X_{<n}, \, X^{<n}\in \D_{<n}$,
$X^n\in \D^n$.
 Then we have the following canonical isomorphisms
\begin{equation}\label{PilPir}
\begin{array}{ll}
\Il(X)\cong X_{<n};\ \
\Ir(X)\cong X^{<n};\\
\Pil\circ \Pi (X)\cong X_n;\ \ \Pir\circ \Pi (X)\cong X^n.
\end{array}
\ \ \epf
\end{equation}


\proofnd of Proposition \ref{tstr}.
 To prove (a) it sufficies to construct a $t$-structure
on $\D$ satisfying  \eqref{defDbo}, \eqref{defDme}: then for another
 $t$-structure $(\D^{\geq 0}_1, \D^{<0}_1)$  such that
$\del^i\in \D_1^{\leq 0}$,
$\nab_i\in \D_1^{\geq 0}$ we have $\D^{< 0}\subseteq \D^{< 0}_1$,
$\D^{\geq 0}\subseteq \D^{\geq 0}_1$, which implies
 $\D^{< 0}= \D^{< 0}_1$, $\D^{\geq 0}= \D^{\geq 0}_1$ (recall
that for a triangulated category $\D$ with a $t$-structure
 $(\D^{\geq 0}, \D^{<0})$
for an object
$X\in \D$ we have $X\in \D^{\geq 0}\iff \Hom(Y,X)=0 \ \ \forall
Y\in \D^{<0}$;  $X\in \D^{< 0}\iff \Hom(X,Y)=0 \ \ \forall
Y\in \D^{\geq 0}$).

 We construct by induction a $t$-structure
on $\D_{\leq n}$ with
\begin{equation}\label{defDbo_upton}
  \D^{\geq 0}_{\leq n}=\< \{\nab^i[d]\ ,\ i\leq n, \; d\leq 0\} \>;
\end{equation}
\begin{equation}\label{defDme_upton}
\D^{< 0}_{\leq n}=\< \{\del_i[d]\ ,\ i\leq n, \; d> 0\} \>;
\end{equation}
since $\D_{<n}=\cupl_{i<n} \D_{\leq i}$ we can assume that
a $t$-structure on $\D_{<n}$ is already defined.

Lemma \ref{2lem} implies that the functors $\D_{<n}\overset{\I}{\To}
\D_{\leq n} \overset{\Pi}{\To} \D_{\leq n}/\D_{<n}$ satisfy the requirements
of \cite{BBD} 1.4.3 (to pass from our notations to those of
\cite{BBD} one should set $\Pi=j^*=j^!$, $\Pil=j_!$, $\Pir=j_*$;
$\I=i_*=i_!$; $\Il= i^*$; $\Ir=i^!$).
 Thus the construction of gluing of $t$-structures
({\it loc. cit.} Theorem 1.4.10) is applicable.

We endow $\D_{<n}$
with the $t$-structure obtained by the induction assumption; and
 $\D_{\leq n}/\D_{<n}\cong \D^n$ with the unique $t$-structure
which has
$\sq_n$ in its core (see Lemma \ref{imposic}).
 Then \cite{BBD} Theorem 1.4.10 provides
$\D_{\leq n}$ with a $t$-structure given by
$$\D_{\leq n}^{\leq 0}=\{X\in \D_{\leq n}\ |\ \Il(X)\in \D_{<n}^{\leq 0}
\ \ \&\ \    \Pi(X)\in (\D_{\leq n}/\D_{<n})^{\leq 0}\};$$
$$\D_{\leq n}^{\geq 0}=\{X\in \D_{\leq n}\ |\ \Ir(X)\in \D_{<n}^{\geq 0}
 \ \ \&\ \     \Pi(X)\in (\D_{\leq n}/\D_{<n})^{\geq 0}\}.$$
In view of \eqref{PilPir}, \eqref{adynelement} we have
$$[\D_{\leq n}^{\leq 0}]=\left[\< \del_n[i]\ | i\geq 0 \>\right]*
[\D_{< n}^{\leq 0}];$$
$$[\D_{\leq n}^{\geq 0}]= [\D_{< n}^{\geq 0}]
* \left[\< \del_n[i]\ | \ i\leq 0 \>\right],$$
which implies \eqref{defDme_upton}, \eqref{defDbo_upton}. The Proposition
is proved. \epf

 We will call the
$t$-structure defined by \eqref{defDbo}, \eqref{defDme} {\it the
$t$-structure of the quasi-exceptional set $\nab$}.

\begin{Rem}
 The $t$-structure of the quasi-exceptional set introduced in Example
\ref{prevrat} is the ``perverse'' $t$-structure \cite{BBD} corresponding to
perversity $p=(p_i)$.
\end{Rem}

\begin{Rem}

 It follows from the axioms of a $t$-structure that the
$t$-structure of an exceptional set $\nab$ can be alternatively
 described as follows.
  For $X\in \D$ we have
\begin{equation}\label{iff}
\begin{array}{ll}
X\in \D^{\geq 0}\iff \Hom^{<0}(\del_i,X)=0\ \forall i\in I;\\
X\in \D^{< 0}\iff \Hom^{\leq 0}(X,\nab^i)=0\ \forall i\in I.
\end{array}
\end{equation}
In the situation of Example \ref{prevrat} \eqref{iff} turns into the usual
definition of a perverse sheaf by a condition on stalks and costalks.
\end{Rem}

We keep the assumptions of Proposition \ref{tstr}. Let $\A$ be the
core of the $t$-structure of the quasi-exceptional set $\nab$;
$\tau$ be the corresponding truncation functors, and
$H^m=\tau_{\leq m}\circ \tau_{\geq m}:\D\to \A$ be the cohomology
functor.

Define $M_i,N^i\in \A$ by $M_i=\tau _{\geq
0}(\del_i)=H^0(\del_i)$, and $N^i=\tau_{\leq
0}(\del^i)=H^0(\del^i)$. Isomorphism \eqref{delnab}
 provides a morphism
$\Phi_i:\del_i\to \nab^i$, which goes to $Id_{\del_i}$ under
\eqref{delnab}, and thus also a morphism $H^0(\Phi_i):M_i\to
N^i$. Also set $\A_{<n}=\A\cap \D_{<n}$, $\A_{\leq n}=\A\cap \D_{\leq n}$;
and let $\A_n$ be the core of the unique $t$-structure on
$\D_{\leq n}/\D_{<n}$ such that $\A_n\owns \sq_n$.

\begin{Prop} \label{irred1}
 Let $L_i$ be the image of $H^0(\Phi_i):M_i\to N^i$. Then $L_i$
is irreducible, and any irreducible object of $\A$ is isomorphic
to $L_i$ for some $i$.

The order on $I$ induces an order on $\{[L_i]\}$, and $\A$ with
this ordering on the set of isomorphism classes of irreducible
objects is a quasi-hereditary category. The canonical morphisms
$\phi_n:M_n\sur L_n$ and $\phi^n: L_n\imbed N^n$ are the standard
and the costandard covers of $L_n$ respectively.
 \end{Prop}

\proof For any $t$-category $\D$ with the core $\A$ and $X\in \A$,
$\del\in \D^{\leq 0}$ we have
\begin{equation}\label{anyt1}
 \Hom(H^0(\del),X)\iso \Hom(\del,X) ;
\ \ \Ext^1_\A(H^0(\del),X)\imbed \Hom(\del,X[1]);
\end{equation}
 and dually for
$\nab \in \D^{\geq 0}$ we have
\begin{equation}\label{anyt2}
\Hom(X,H^0(\nab))\iso \Hom(X,\nab);\
\
 \Ext^1_\A(X,H^0(\nab))\imbed \Hom(X,\nab[1]).
 \end{equation}
  Thus we have
\begin{equation}\label{Mort}
\Hom(M_n,X)=0=\Ext^1(M_n,X);
\end{equation}
\begin{equation}\label{ortN}
\Hom(X,N^n)=0=\Ext^1(X,N^n)
\end{equation}
for  $X\in \A_{<n}$. Let us now show that $L_n$ is simple in $\A$.
Assume that $0\to X\to L_n\to Y\to 0$ is a short exact sequence.
Pick the minimal  $i$ such that $X,Y\in \A_{\leq i}$. If $i>n$ we get a
contradiction because applying the exact functor
$\Pi_i$ to the exact sequence we get
an exact sequence
$$0\to X\mod \D_{<i} \to 0 \to Y\mod \D_{<i}\to 0$$
in $\A_i$, which shows that $X\mod \D_{<i}=0=Y\mod \D_{<i}$,
so $X,Y \in \D_{\leq j}$ for some $j<i$. Thus $i=n$; so we have
an exact sequence in $\A_n$
$$0\to X\mod \D_{<n} \to \sq_n\to Y\mod \D_{<n}\to 0.$$
Since $\sq_n$ is irreducible by Corollary \ref{el1}, we see that
either $X\in \A_{<n}$ or $Y\in \A_{<n}$. However, $X$ is a
subobject of $N^n$, while $Y$ is a factor-object of $M_n$; thus we
get a contradiction with either \eqref{Mort} or \eqref{ortN}.

We claim that
\begin{equation}\label{Lgen}
\A=\< L_n \ |\ n\in I\>.
\end{equation}
Notice that \eqref{Lgen} implies the statement of the Lemma: since
a sequences $0\to X\overset{f}{\to} Y\overset{g}{\to} Z\to 0$ in
$\A$ is exact iff there exists a distinguished triangle
$$X\overset{f}{\to} Y\overset{g}{\to} Z\overset{h}{\to} X[1],$$
\eqref{Lgen} shows that any object of $\A$ has a finite filtration
with every subquotient isomorphic to $L_n$ for some $n$. To check
\eqref{Lgen} observe that the isomorphism $L_n\cong \nab^n\mod
\D_{<n}$ implies by induction that $L_i$, $i\leq n$ generated
$\D_{\leq n}$ as a triangulated category. So \eqref{Lgen} follows
from Lemma \ref{imposic}. \epf

\begin{Rem}\label{irred2}
 Fix $i\in I$, and set,  $(\nab')^j=\nab^j[1]$,
$(\nab'')^j=\nab^j[-1]$ for $j<i$; $(\nab')^j= \nab^j=(\nab'')^j$
for $j\geq i$. Then $\nab'$, $\nab''$ are dualizable
quasi-exceptional sets; let $\tau'$, $\tau''$ be the truncation
functors for the corresponding $t$-structures. One can show that
$$\tau'_{\leq 0}(\nab^i)\cong L_i\cong \tau''_{\geq 0}(\del_i).$$
\end{Rem}

\begin{Rem} In the situation of Example \ref{prevrat}
Proposition \ref{irred1} provides the standard description of a
Goresky-MacPherson IC-sheaf $j_{!*}(\L)$ (where $\L$ is a local
system) as the image of the canonical morphism $H^{p,0}(j_!(\L))
\to H^{p,0}(j_*(\L))$; while Remark \ref{irred2} describes $j_{!*}(\L)$ as a
result of successive applications of the direct image and
truncation functors, cf \cite{BBD}, Proposition 2.1.11 (cf also {\it loc.
cit.} 2.1.9).
\end{Rem}

\section{Main result}
The pull-back and push-forward functors for coherent
 sheaves are understood to be the corresponding
 derived functors, unless stated otherwise.

We return to the set-up and notations of the introduction. In particular
$\pi:\Ntil=T^*(G/B) \to \N$ is the moment map from the cotangent bundle
of the flag variety  $G/B$ to the nil-cone
(the Springer-Grothendieck resolution);
 let also $p:\Ntil\to G/B$ be the
projection.

 From now on we set $\D=\Db(Coh^G(\N))$.

 For a weight
$\lambda\in \Lambda$ let $\O_{G/B}(\lambda)$
 be the corresponding  $G$-equivariant line bundle on $G/B$
 (thus $\O_{G/B}(\lambda+\rho)$ is ample
for $\lambda\in \Lambda^+$); for a parabolic $P\subset G$ we
will  write $\O_{G/P}(\lambda)$ for the unique equivariant line
bundle on $G/P$ whose pull-back to $G/B$ is $\O_{G/B}(\lambda)$ if
such a line bundle exists; for a variety $X$ with a map $f:X\to
G/P$ we will denote $f^*(\O_{G/P}(\lambda)) $ by $\O_X(\lambda)$;
we will write $\F(\lambda)$ instead of $\F\otimes \O_X(\lambda)$
for $\F\in \Db(Coh_X)$ etc. For $\lambda \in \Lambda^+$ we set
$V_\lambda=H^0(G/B,\O_{G/B}(\lambda))$.

 We define $A_\lambda\in \D $ by
$A_\lambda= R\pi_* (
 \O_{\Ntil}(\lambda) )$.

Let $W$ be the Weyl group.
 For $\lambda\in \Lambda$ we denote
its $W$ orbit by $W(\lambda)$; let $conv(\lambda)$ be the intersection
of the convex hull of $W(\lambda)$ with $\Lambda$, and $conv^0(\lambda)$
be the complement to $W(\lambda)$ in $conv(\lambda)$.

For a subset $S\subset \Lambda$ let
$$\D_S=\D_{\{A_\lambda, \ \lambda\in S\}}
=\< A_\lambda[n],\ \lambda\in S,\ n\in \Zet
\>$$
 be the  triangulated
subcategory of $\D$ generated by $A_\lambda$, $\lambda\in S$.

\begin{Prop}\label{W}
For $w\in W$ there exists a canonical isomorphism
\begin{equation}\label{isow}
A_\lambda\cong A_{w(\lambda)} \mod \D_{conv^0(\lambda)}.
\end{equation}
\end{Prop}

 \proofnd is a variation of a classical argument going back 
at least 
to \cite{De}.

Let $\alpha$ be a simple root, $s_\alpha\in W$ the corresponding
simple reflection. It sufficies to construct \eqref{isow} for $w=s_\alpha$
and  $\lambda\in \Lambda$  such that
$s_\alpha(\lambda)=\lambda - n \alpha$, $n>0$.

Let $\pr_\alpha:G/B\to G/P_{\alpha}$ be the projection, where
$P_\alpha$ is the minimal parabolic corresponding to $\alpha$. Let
$G'\to G$ be the universal covering, and let $\Lambda'\supset
\Lambda$ be the weight lattice of $G'$. There exists $\lambda'\in
\Lambda'$ such that $s_\alpha(\lambda')=\lambda'-(n-1)\alpha$.
Set  $\Ve_{\lambda'}=\pr_{\alpha}^*\pr_{\alpha
*}\left(\O_{G/B}(\lambda')\right)$. Thus $\Ve_{\lambda'}$ is
$G'$-equivariant  vector bundle; it has a $G'$-invariant
filtration with subquotients $s_\alpha(\lambda')$,
$s_\alpha(\lambda')+\alpha,\dots, \lambda'$. We claim that
\begin{equation}\label{veve}
\pi_*\left(  p^*(\Ve_{\lambda'})(\lambda - \lambda')\right)
\cong \pi_* \left(  p^*(\Ve_{\lambda'})\left(
\lambda - \lambda'-\alpha \right)\right).
\end{equation}
Since $p^*(\Ve_{\lambda'})(\lambda - \lambda')$ is a
$G$-equivariant vector bundle on $\Ntil$ equipped with a filtration
whose subquotients are $\O_\Ntil(\lambda -k\alpha)$, $k =0,\dots
,n-1$, we have $$[\pi_*\left(  p^*(\Ve_{\lambda'})(\lambda -
\lambda')\right)]\in [A_{\lambda -(n-1)\alpha}]*\cdots *
[A_\lambda],$$ and $$\pi_*\left(  p^*(\Ve_{\lambda'})(\lambda -
\lambda')\right)\cong A_\lambda \mod \D_{conv^0(\lambda)}.$$
Similarly $$[\pi_* \left(
p^*(\Ve_{\lambda'})\left(\lambda -
\lambda'-\alpha \right)\right)]\in [A_{\lambda-n\alpha }]*\cdots
*[A_{\lambda -\alpha}],$$ hence $$\pi_* \left(
p^*(\Ve_{\lambda'})\left(s_\alpha(\lambda -
\lambda')\right)\right)\cong A_{\lambda-n\alpha}=
A_{s_\alpha(\lambda)}\mod
\D_{conv^0(\lambda)}.$$ Thus \eqref{veve} yields \eqref{isow}.

It remains to  check \eqref{veve}.
 Set $\Ntil_\alpha=T^*(G/P_\alpha)\times_{G/P_{\alpha}}
G/B$; the differential of $\pr_\alpha$
provides a closed imbedding $           
 \Ntil_\alpha \imbed
\Ntil$. We have an exact
sequence in $Coh^G(\Ntil)$
\begin{equation}\label{ese}
0\to \O_{\Ntil}(\alpha)\to \O_{\Ntil} \to \O_{\Ntil_\alpha} \to 0.
\end{equation}
Tensoring it with $p^*(\Ve_{\lambda'})(\lambda - \lambda'-\alpha)$
we see that to check \eqref{veve} it suffices to verify that
$$\pi_*\left( \O_{\Ntil_\alpha}\otimes p^*(\Ve_{\lambda'})(\lambda -
\lambda'-\alpha) \right)=0.$$ We claim that in fact a stronger equality
$$\pi_*'\left( \O_{\Ntil_\alpha}\otimes p^*(\Ve_{\lambda'})(\lambda -
\lambda'-\alpha) \right)=0$$ holds, where $\pi'$ is the projection
$$\pi':\Ntil_\alpha = T^*(G/P_\alpha)\times_{G/P_\alpha} G/B \To
T^*(G/P_\alpha).$$ Indeed, the fibers of $\pi'$ are projective
lines, and $\O_{\Ntil_\alpha}\otimes p^*(\Ve_{\lambda'})(\lambda -
\lambda'-\alpha)$ is readily seen to be isomorphic to a sum of
several copies of $\O_{\mathbb P^1} (-1)$ when restricted to any
fiber of $\pi'$. \epf

\begin{Prop}\label{Agood}
a) Let $\lambda\in\Lambda^+$. Then $\Hom(A_\lambda,A_\lambda)=\Ce$, and
$\Hom^{<0}(A_\lambda,A_\lambda)=0$.
Also for any $\mu\in \Lambda$ we have:
\begin{equation}\label{HomAA}
\Hom^\bu(A_\mu, A_\lambda)=0 {\ \ \rm if \ \ } \lambda\not \in conv(\mu).
\end{equation}

b) For $\lambda, \mu \in \Lambda^+$,  $\lambda\ne \mu$
we have
\begin{equation}\label{HomAwA}
\Hom^\bu(A_{w_0(\mu)},A_\lambda)=0,
\end{equation}
where $w_0\in W$ be the element of maximal length.

\end{Prop}

 We will need the following known fact.

\begin{Fact}\label{fact} \cite{Br}, \cite{K}
a) For dominant $\lambda$ we have $H^i(\Nt, \O(\lambda))=0$ for $i\ne 0$,
and
\begin{equation}\label{dimHom}
\dim \Hom_G(V_\mu, H^0(\Nt, \O(\lambda))=n_\lambda^\mu,
\end{equation}
where $n_\lambda^\mu$ is the multiplicity of weight $\lambda$ in
$V_\mu$.

b) 
 $R^\bu\pi_*\O_\Nt=\O_\N$. \epf
\end{Fact}

\begin{Rem}
We will only use \eqref{dimHom} in the case when
 $\mu \not \in conv(\lambda)$, so both sides vanish.
\end{Rem}

\proofnd of Proposition \ref{Agood}(a).
By Fact \ref{fact} $A_\lambda$ is a sheaf (rather than a complex) for
$\lambda\in \Lambda^+$; thus of course $\Hom^{<0}(A_\lambda, A_\lambda)=0$.
Also $A_\lambda$ is torsion free and has generic rank 1,
hence $\Hom^G(A_\lambda,A_\lambda)=\Ce$, because $\N$ has an open orbit.
It remains to check \eqref{HomAA}.

From  Fact  \ref{fact} it follows that
if $\lambda\in \Lambda^+$, then  $R^i\Gamma(A_\lambda)=
H^i(\Nt, \O(\lambda))=0$ for $i\ne 0$, and $\Hom_G(V_\mu, \Gamma(A_\lambda))
=\Hom_G(V_\mu, \Gamma(\Nt, \O(\lambda) ))=0$ unless $\lambda\in conv(\mu)$.
Thus we get
\begin{equation}\label{HomV}
\Hom^\bu (V_\mu\otimes \O, A_\lambda)
=\Hom_G(V_\mu, R^\bu\Gamma(A_\lambda))
=0\ \ \ {\rm if}\ \ \lambda\not \in conv(\mu).
\end{equation}
Introduce a (nonstandard) order on $\Lambda$ by
$\nu_1\preceq \nu_2$ if $\nu_1\in conv(\nu_2)$.
We fix $\lambda$, and proceed by induction
in $\mu$ in that order. We can assume that \eqref{HomAA} holds
for all $\mu'\in conv^0(\mu)$. Now notice that $V_\mu\otimes \O_{G/B}$
carries a filtration whose associated graded
is  $$
\gr(V_\mu\otimes \O_{G/B})=\oplusl\O_{G/B}(\nu)^{\oplus n_\nu^\mu}.$$
 Hence
$$
[V_\mu\otimes \O_\N]=[\pi_*p^*(V_\mu\otimes \O_{G/B})]
\in \{ [A_{\nu_1}^{n_{\nu_1}^\mu} ] \} * \cdots
 \{ [A_{\nu_k}^{n_{\nu_k}^\mu }] \},
$$
where $\nu_1,\dots,\nu_k$ are weights on $V_\mu$.
The induction assumption says that $\Hom^\bu (A_\nu,A_\lambda)=0$
for all $\nu\in {conv^0(\mu)}$. Thus the last equality
implies that
\begin{multline}\label{twinkletwinklelittlestar}
[0]=[\Hom^\bu(V_\mu\otimes \O_\N, A_\lambda)]\in
\{[\Hom^\bu(A_{\nu_k}^{ n_{\nu_k}^\mu}, A_\lambda)]\}*\cdots *
\{[\Hom^\bu(A_{\nu_k}^{ n_{\nu_1}^\mu}, A_\lambda)]\} =\\
\{[\Hom^\bu(A_\mu, A_\lambda)]\}* \{[\Hom^\bu(A_{w_1(\mu)},
A_\lambda)]\}*\cdots *\{[\Hom^\bu(A_{w_s(\mu)}, A_\lambda)]\},
\end{multline}
where  $\mu, \, w_1(\mu),\dots , w_s(\mu)$ are extremal weights of
$V_\mu$ (here we view $\Hom^\bu$ as an object of
$D^+(Vect_{\Ce})$). By Proposition \ref{W}, $A_{w(\mu)}\cong A_\mu \mod
\D_{conv^0(\mu)}$, thus by the induction assumption,
$\Hom^\bu(A_{w(\mu)}, A_\lambda)\cong \Hom^\bu(A_\mu, A_\lambda)$
for all $w\in W$. So \eqref{twinkletwinklelittlestar} can be
rewritten as
$$[0]\in \{[\Hom^\bu(A_\mu, A_\lambda)]\}*\{[\Hom^\bu(A_\mu, A_\lambda)]\}*
\cdots * \{[\Hom^\bu(A_\mu, A_\lambda)]\},$$
where the number of terms in the right-hand side is the number
of extremal weights in $V_\mu$.
Now \eqref{HomAA}
follows from the next standard Lemma, applied to $\A=Vect$,
$V=\Hom^\bu(A_\mu,A_\lambda)$. \epf

\begin{Lem}\label{notcancel} Let $V\in D^+ (\A)$ for an abelian category $\A$.
If $$[0]\in \{[V]\}*\{[V]\}*\cdots *\{[V]\}$$ (where $V$ is
repeated $n$ times, $n\geq 1$), then $V=0$.
\end{Lem}

\proof Otherwise, if $i$ is  minimal, such that $H^i(V)\ne 0$,
then $H^i(V)\imbed H^i(X)$ for $[X]\in \{[V]\}*\{[V]\}*\cdots
*\{[V]\}$. \epf

\medskip

\proofnd of Proposition \ref{Agood}(b).
 If $\lambda\not \in conv(\mu)$, then \eqref{HomAwA} follows from
\eqref{HomAA}. Otherwise, $w_0(-\mu)\not \in conv(-\lambda)$.
Recall that the Grothendieck-Serre duality $\GS$
is an anti-autoequivalence of $\D$, such that
$\GS(A_\lambda)=A_{-\lambda}[\dim \N]$, see section
\ref{32} below.
Thus we have $$\Hom^\bu(A_{w_0(\mu)},A_\lambda)=\Hom^\bu(\GS(A_\lambda),
\GS(A_{w_0(\mu)}))=\Hom^\bu(A_{-\lambda},A_{w_0(-\mu)}),$$
which again vanishes by \eqref{HomAA}. \epf

\begin{Prop}\label{Dgen}
 $\D$ is generated by $A_\lambda$, $\lambda\in \Lambda^+$.
\end{Prop}

\begin{Lem}\label{absgen}\footnote{Our assumptions in this Lemma coincide with
those of \cite{CG}; the statement is apparently true in a more general
situation.}
Let $X$ be an algebraic variety over $\Ce$, and $p:Y\to X$ be a vector
bundle; let $G$ be a linear algebraic group acting on $X$, $Y$,
so that $p$ is $G$-equivariant. Then $D^b(Coh^G(Y))$ is generated
as a triangulated category by objects of the form $p^*(\F)$, $\F\in Coh^G(X)$.
\end{Lem}

\proof See e.g. \cite{CG}, p. 266 (last paragraph). \epf

\begin{Cor}\label{Ntilgen}
 $D^b(Coh^G(\Ntil))$ is generated by the objects
$\O_{\Ntil}(\lambda)$, $\lambda \in \Lambda$.
\end{Cor}

\proof
The category $Coh^G(G/B)$ is identified with the category of
representations of $B$; in particular, any object of $Coh^G(G/B)$ is a vector
bundle, which carries a filtration with subquotients being $\O_{G/B}(\lambda)$,
$\lambda\in \Lambda$. Now apply Lemma \ref{absgen} to $X=G/B$, $Y=\Ntil$. \epf

\begin{Lem}\label{pisur}
The image of the functor $\pi_*:\Db(Coh^G(\Ntil))\to \D$ 
generates $\D$ as a triangulated category.
\end{Lem}
\proof (cf e.g.  \cite{O}, Lemma 2.2).
It is enough to show that for $\F\in \D$ there exists
$\Ftil\in D^b(Coh^G(\Ntil))$ and a morphism $\phi:\F\to \pi_* (\Ftil)$
such that the support of its cone $Cone(\phi)$ is strictly smaller than
the support of $\F$.
We can assume that $\F\in Coh^G(\N)\subset \D$, and also that the
scheme-theoretic support of $\F$ is reduced.
 Let
$O\subset \N$ be a $G$-orbit which is open in the
support of $\F$, and $\Obar$ be its closure.
 It is well known  
that there exists a
 parabolic subgroup $P_O$, and
  a $G$-equivariant  subbundle $\Ntil_O\subset T^*(G/P_O)$, such that
  $\pi_O:{\Ntil_O}\to \Obar$ is birational (and thus is a resolution of
  singularities of $\Obar$); here $\pi_O$ is the restriction to $\Ntil_O$
of the moment map $T^*(G/P_O)\to \N$ (more presicely, for
$x\in O$ one can define $P_O$, $\Ntil_O$ by  ${\mathfrak p}_O=\g_{\geq 0}$,
$\Ntil_O=\g_{\geq 2}\times_{P_O} G$; here 
$\g_{\geq i}$ are the terms of the Jacobson-Morozov-Deligne filtration
on the Lie algebra $\g$ associated to the nilpotent operator ${\rm ad}(x)$,
 see \cite{Weil2}, 1.6, and ${\mathfrak p}_O$ is the Lie algebra of $P_O$).
 Let $\itil$ be the imbedding $\Ntil_O\times_{G/P_O}G/B
\imbed \Ntil$, and set $\Ftil=(\pi\circ \itil)^*(\F)$ (the non-derived
pull-back). We have a canonical adjunction morphism
$\F\to (\pi\circ \itil)_*(\Ftil)$ (where we again consider the non-derived
direct image). The composition
$$\F\to (\pi\circ \itil)_*(\Ftil)=R^0(\pi\circ \itil)_*(\Ftil) \to
R\pi_*(\itil_* \Ftil)$$ is an isomorphism on $O$, because the
fiber of $\pi\circ \itil$ over a point of $O$ is $P_O/B$,
(the flag variety for the Levi subgroup),
and the structure sheaf of $P_O/B$ is acyclic; hence the cone of
this composition has smaller support. \epf

\proofnd of Proposition \ref{Dgen}. It  follows directly from Lemma
\ref{pisur} and Corollary \ref{Ntilgen} that $\D$ is generated by
$A_\lambda$, $\lambda\in \Lambda$. So it is enough to show that
for  $\lambda\in \Lambda^+$ the category $\D_{conv(\lambda)}$ is
generated by $A_\mu$, $\mu \in \Lambda^+\cap conv(\lambda)$. This
follows by induction in $\lambda$ (with respect to  the standard
 partial order on $\Lambda^+$) from Proposition \ref{W}. \epf

 Propositions
\ref{Agood},   \ref{Dgen} and \ref{W} yield the following
\begin{Thm}\label{Aqe}
Let us equip $\Lambda^+$ with any total ordering $\leq$
compatible with the standard partial order (i.e $\lambda \in
conv(\mu)\Rightarrow \lambda\leq \mu$). Then the set
$\{A_{\lambda}\ |\ \lambda \in \Lambda^+\}$ is a quasi-exceptional
set generating $\D$. The set $$\{A_{w_0(\lambda)}\ |\ \lambda \in
\Lambda^+\}$$ is a dual quasi-exceptional set.
\end{Thm}

\begin{Rem}
The set $\{A_\lambda\}$ is not exceptional. For example, one can
show that if $G=SL(2)$, and $\lambda=1$, then $\Ext^i(A_\lambda,
A_\lambda)\cong \Ce$ for all $i\geq 0$. It is also easy to see,
that $\D=\Db(Coh^G(\N))$ is not generated by any exceptional set
(for otherwise $\Hom^\bu(X,Y)$ would be finite dimensional for all
$X,Y\in \D$, while this is not so in the above example
$X=Y=A_1$,  $G=SL(2))$. Notice, however, that the ``larger''
category $\Db(Coh^G(\Nt))$ is generated by the set
$\O(\lambda)$, $\lambda\in \Lambda$, which can be shown
to be exceptional for any  ordering on
$\Lambda$, which is compatible with the standard partial order.
\end{Rem}

\begin{Rem}\label{deltil}
Let $\Pi_\lambda:\D_{\leq \lambda}\to \D_{\leq \lambda}/\D_{<\lambda}$ be the
projection; by Theorem \ref{Aqe}, Lemma \ref{2lem} (d) we have left and right
adjoint functors $\Pil_\lambda:\D_{\leq \lambda}/\D_{<\lambda}\to \D_
{A_{w_0(\lambda)}}$, and 
 $\Pir_\lambda:\D_{\leq \lambda}/\D_{<\lambda}\to \D_
{A_{\lambda}}$. For $\lambda\in \Lambda^+$ set 
$$\begin{array}{ll}
\bA_\lambda= \Pir_\lambda(V_\lambda\otimes \O_\N);\\
\bA_{w_0(\lambda)}= \Pil_\lambda(V_\lambda\otimes \O_\N).
\end{array}$$
Then 
we have
$$
\Hom^\bu(\bA_{w_0(\lambda)}, A_\mu)= 0 = \Hom^\bu(A_{w_0(\lambda)}, \bA_\mu)
$$
for $\lambda\ne \mu \in \Lambda^+$, and
$$\Hom^\bu(\bA_{w_0(\lambda)}, A_\lambda)=\Ce =
 \Hom^\bu(A_{w_0(\lambda)}, \bA_\lambda),$$
where the latter equality follows from \eqref{dimHom}.

We claim that $\bA_\lambda$, $\bA_{w_0(\lambda)}$ admit the following
geometric description. 

For $\lambda\in \Lambda$ let $P_\lambda$ be the largest parabolic such that
$\O_{G/B}(\lambda)$ is isomorphic to the pull-back of a line bundle on 
$G/P_\lambda$; let $\p_\lambda$ be its Lie algebra.
Set $\gtil_\lambda = \p_\lambda \times_{P_\lambda} G$, and
let $\pi_\lambda:\gtil_\lambda\to \g$ be the projection.
 Then we have
\begin{equation}\label{geomop1}
\bA_\lambda \cong i^*\pi_{\lambda*}(\O_{\gtil_{\lambda}}(\lambda));  
\end{equation}
\begin{equation}\label{geomop2}
\bA_{w_0(\lambda)} \cong i^*\pi_{\lambda*}\left(\O_{\gtil_{\lambda}}
\left(w_0    (\lambda)\right)\right),
\end{equation}
where $i:\N\to \g$ is the imbedding, and
$\O_{\gtil_{\lambda}}(\lambda)$ is defined by means of the obvious projection
$\gtil_\lambda\to G/P_\lambda$. 

Indeed, the familiar morphism $V_\lambda\otimes\O_{G/P_\lambda}\to 
\O_{G/P_\lambda}(\lambda)$ yields a morphism $V_\lambda \otimes 
\O_{\gtil_\lambda}\to \O_{\gtil_\lambda}(\lambda)$, and hence 
also morphisms
 $$V_\lambda\otimes \O_\g \to \pi_{\lambda*}(\O_{\gtil_{\lambda}}(\lambda)),$$
 $$V_\lambda\otimes \O_\N
\to i^*\pi_{\lambda*}(\O_{\gtil_{\lambda}}(\lambda));$$
and thus a morphism $\phi$ from the left-hand side to the right-hand side of
\eqref{geomop1}. Since both objects in question lie in $\< A_\lambda\>$
and have length $\dim (H^\bu(G/P_\lambda))$, it sufficies to check
that this morphism is injective. This would follow if we show that
the composition $A_{w_0(\lambda)}\to V_\lambda\otimes \O_\N\overset{\phi}{\To}
 \bA_\lambda$
is non-zero, where the first arrow is the only (up to a constant) non-zero 
morphism $A_{w_0(\lambda)}\to V_\lambda$. Thus it sufficies to see that
$\phi|_{\N_0}$ is surjective, where $\N_0\subset \N$ is the open orbit.
Surjectivity of $\phi|_{\N_0}$ follows from the next Lemma, which
is an unpublished result of Bram Broer.
Finally
\eqref{geomop2} follows from \eqref{geomop1} by Grothendieck-Serre duality. 
\end{Rem}

\begin{Lem}
Let $e\in \g$ be a regular nilpotent, let $\gtil_\lambda^e$ be the preimage
of $e$ under $\pi_\lambda$, and $(G/P_\lambda)^e$ be the
image of $\gtil_\lambda^e$ in $G/P_\lambda$ (thus  $(G/P_\lambda)^e$
is a nilpotent scheme of length $\dim H^*(G/P_\lambda)$).

Then the restriction map 
\begin{equation}\label{restri}
V_\lambda=
\Gamma(G/P_\lambda, \O(\lambda))
\to \Gamma((G/P_\lambda)^e, \O(\lambda))
\end{equation}
is surjective. \epf
\end{Lem}

\begin{Rem}
Victor Ginzburg pointed out to us that the surjection \eqref{restri}
 probably admits the following
alternative description. One can realize $V_\lambda$ as the total cohomology
of an irreducible perverse sheaf $IC_\lambda$
on the affine Grassmanian $\Gr$ of the Langlands dual group
$^LG$, equivariant under
the maximal bounded subgroup $^LG(O)$ in the loop group $^LG(K)$,
see \cite{Gi}, \cite{MV}. Ginzburg
conjectures that one can identify $ \Gamma((G/P_\lambda)^e, \O(\lambda))$
with cohomology (with constant coefficients) of the open $^LG(O)$ orbit
$\Gr_\lambda$
in the support of $IC_\lambda$, so that \eqref{restri} is identified with
the restriction map
\begin{equation}\label{restri1}
H^\bu(IC_\lambda) \to H^\bu(IC_\lambda|_{\Gr_\lambda})= H^\bu (\underline{\Ce}
_{\Gr_\lambda}[\dim (\Gr_\lambda)]),
\end{equation}
where $\underline{\Ce}
_{\Gr_\lambda}$ is the constant sheaf on $\Gr_\lambda$.

Notice that it is easy to see that  $$H^\bu(\Gr_\lambda)
\cong H^\bu(^LG/^LP_{\lambda\check{\ }})\cong 
 H^\bu (G/P_\lambda)\cong \Gamma((G/P_\lambda)^e, \O),$$
where $\lambda\check{\ }$ is a weight of $^LG$ obtained from $\lambda$
by means of an invariant quadratic form on $\g$.
Thus at least the dimensions of the target spaces in 
\eqref{restri} and \eqref{restri1} coincide. 
\end{Rem}

\subsection{Comparison with perverse coherent $t$-structure}\label{32}
Recall the coherent perverse $t$-structure on $\D$, corresponding
to the middle perversity, $p(O)=-\frac{\dim(O)}{2}$ for a $G$-orbit $O\subset
\N$, see \cite{izvrat}.
 We let $\D^{p,>0}$, $\D^{p,\leq 0}$ denote the corresponding
positive and negative subcategories, and $\P=\D^{p,\geq 0}\cap \D^{p,\leq 0}$
be its core.

Let $\GS:\D\to \D^{op}$ be the Grothendieck-Serre duality;
$\GS:X\mapsto \underline {RHom} (X, \DC)$, where $\DC$
 is the {\it equivariant dualizing complex}, cf. \cite{izvrat}, Definition 1
(we assume that the dualizing complex is normalized in the standard
way, i.e. $\DC=pr^!(\underline\Ce)$, where $pr$ is the projection to
$Spec(\Ce)$).

Set $d=\dim(\N)/2$.
\begin{Lem}\label{DinD}
We have $A_\lambda [d]
\in \P$ for all $\lambda$.

\end{Lem}

\proof We have $\GS(A_\lambda)[d]=A_{-\lambda}[d]$, because
duality commutes with proper direct images, and the (equivariant)
dualizing
sheaf on $\Ntil$ is isomorphic to $\O_\Ntil[2d]$. Thus it is
enough to check conditions on stalks on $A_\lambda[d]$; i.e.
for an orbit $O\subset \N$ we have to see that
$$i_{O_{\rm gen}}^*(A_{-\lambda})\in \D^{\leq p(O)+d}(\O_{O_{\rm gen}}-mod),$$
where $i_{O_{\rm gen}}:O_{\rm gen}\imbed O$
is the imbedding of the generic point of $O$ into $\N$
(see \cite{izvrat}, Definition 2, Lemma 5(a)).
This follows from the two well-known facts:
that $\pi$ is a semi-small morphism,
(i.e. $\dim (\Ntil\times_\N\Ntil)=\dim \Ntil$,
so that $\dim(\pi^{-1}(x))
\leq \frac{1}{2}\codim(O)=p(O)+d$
for an orbit $O\subset \N$, $x\in O$); and
 that  the homological dimension of the direct image functor
$\pi_*$ for coherent sheaves under a proper morphism $\pi$
of algebraic varieties over a field equals the dimension of
$\pi$ (maximal dimension of a fiber of $\pi$)
(see e.g. \cite{CG}, 3.3.20, 8.9.19; and
\cite{Ha} Corollary 11.2 respectively). \epf

For $(O,L)\in \nildat$  (notations of the Introduction) let
$\L$ be the $G$-equivariant vector bundle on $O$ corresponing to $L$,
and $j:O\imbed \N$ be the imbedding. We set
 $$IC_{O,L}=j_{!*}(\L[p(O)])\in D^b(Coh),$$
 see \cite{izvrat}, 3.2 for the definition of the minimal (Goresky-MacPherson) 
extension functor $j_{!*}$ for coherent sheaves.

\begin{Rem} We do not know an explicit description of the object
$IC_{O,L}$ in general; however, in the particular case when $L=\Ce$
is the trivial representation they are easy to describe. Namely, we
claim that for any orbit $O\overset{j}{\imbed} \N$ we have
\begin{equation}\label{IC}
IC_{O,L}\cong j_*\O_O[p(O)]=N_*(\O  
)[p(O)],
\end{equation} where $j_*$ stands for the {\it non-derived} direct image,
and $N$ is the normalization morphism for $\Obar$ 
(cf. Conjecture 4 in \cite{O}). 
 Indeed, the result of \cite{Hi}, \cite{Pa} implies that the normalization
of $\Obar$ is Cohen-Macaulay; hence $$\GS(N_*(\O)[p(O)])\in Coh^G(\N)[p(O)],$$
which yields \eqref{IC}.

Ostrik pointed out to us that a similar statement is probably true for
$(O,L)\in \nildat$ if $L$ has finite image, due to the result of \cite{Br1}. 
\end{Rem}

\begin{Cor}\label{DisD}
a) The perverse $t$-structure  on $D^b(Coh^G(\N))$
 corresponding to
the middle perversity (\cite{izvrat}, Theorem 1; Example 1)
coincides with the $t$-structure of the
dualizable quasi-exceptional set $\nab_\lambda=A_{\lambda}[d]$.

b) The core $\P$ of this $t$-structure is a quasi-hereditary category.

The set of isomorphism classes of irreducible objects in $\P$ equals
$\{ [IC_{O,L}]\ \ |\ \ (O,L)\in \nildat\}$.

The set of isomorphism classes of costandard objects equals 
$\{[\nab^\lambda]=
[A_\lambda[d]]\ |\ \lambda\in \Lambda^+\}$; and that of standard objects
equals $\{[\del_\lambda]=[A_{w_0(\lambda)}[d]]\ |\ \lambda\in \Lambda^+\}$.
\end{Cor}
\proof (a) follows directly from Lemma \ref{DinD}, and Proposition
\ref{tstr}.
First statement in part (b) is a particular case of Corollary 4 in \cite{B}.
The rest follows from Proposition \ref{irred1}. 
\epf

\begin{Cor} a) The Grothendieck group $K^0(\D)$ is a free abelian
group; each of the sets $\{[A_\lambda] \ |\ \lambda\in
\Lambda^+\}$, $\{[IC_{O,L}]\ |\ (O,L)\in \nildat\}$ forms a basis
of this group.

b) There exists a unique  bijection between $\Lambda^+$ and
$\nildat$, $\lambda\mapsto (O_\lambda, L_\lambda)$ satisfying
either of the following equivalent properties.

(i) $\Hom(IC_{O_\lambda, L_\lambda}, A_\lambda[d])\ne 0$.

$(i')$  $\Hom(A_{w_0(\lambda)}[d],IC_{O_\lambda, L_\lambda})\ne 0$.

(ii) There exists a morphism $IC_{O_\lambda, L_\lambda}\to A_\lambda[d]$
whose cone lies in $\D_{conv^0(\lambda)}$.

$(ii')$ There exists a morphism $A_{w_0(\lambda)}[d]\to
IC_{O_\lambda, L_\lambda}$
whose cone lies in $\D_{conv^0(\lambda)}$.

(iii) $IC_{O_\lambda, L_\lambda}\in \D_{conv(\lambda)}$ for all
$\lambda$.

(iv) $[IC_{O_\lambda, L_\lambda}]$ lies in the span of $[A_\mu]$,
$\mu \in \Lambda^+\cap conv(\lambda)$, i.e. the transformation
matrix between the two bases is upper triangular.

(v) $[IC_{O_\lambda, L_\lambda}]-(-1)^d [A_\lambda]$ lies in
the span of $[A_\mu]$, $\mu \in conv^0(\lambda)\cap \Lambda^+$.
\epf
\end{Cor}

\end{document}